\documentclass{article}
\usepackage{amsfonts}
\usepackage{amssymb}
\usepackage{graphicx}
\usepackage{amsmath}

\newtheorem{theorem}{Theorem}

\newtheorem{corollary}[theorem]{Corollary}

\newtheorem{example}[theorem]{Example}

\newtheorem{lemma}[theorem]{Lemma}

\newtheorem{remark}[theorem]{Remark}

\newenvironment{proof}[1][Proof]{\textbf{#1.} }{\ \rule{0.5em}{0.5em}}

\input{tcilatex}

\begin{document}

\title{Weak Convergence of Laws on $\mathbb{R}^{K}$ with Common Marginals }
\author{Alessio Sancetta\thanks{\textbf{Aknowledgments: }Partially supported
by the ESRC Award RES 000-23-0400. \textbf{Address for Correspondence:}
Faculty of Economics, Sidgwick Avenue, University of Cambridge, Cambridge
CB3 9DE, UK. Tel.: +44-(0)1223-335272; E-mail:
alessio.sancetta@econ.cam.ac.uk. } \\
University of Cambridge}
\maketitle

\begin{abstract}
We present a result on topologically equivalent integral metrics (Rachev,
1991, M\"{u}ller, 1997) that metrize weak convergence of laws with common
marginals. This result is relevant for applications, as shown in a few
simple examples.


\textbf{Classification AMS:} 60B10, 60E05.
\end{abstract}

\section{Introduction}

In applications, it is often necessary to resort to approximations in the
study of the properties of stochastic models. The original process can be
replaced by a simpler one whose characteristics are already known or easier
to study. This requires some stability of the model, which is usually
represented in terms of integral probability metrics (Rachev, 1991, M\"{u}%
ller, 1997). Among the many possible situations, a typical one is the
difficulty to deal with the dependence properties of stochastic processes.
For example, the study of the central limit problem for the standardised
partial sum of dependent random variables requires the introduction of
dependence conditions (mixing conditions, e.g. Doukhan, 1994, or weak
dependence conditions, e.g. Doukhan and Louhichi, 1999). These dependence
conditions are used to relate past and future realisations of some
stochastic process with some "independent version": an example in the text,
below, will make this statement less elusive. In this case as in many other
applications, the univariate marginal distributions of the approximation are
the same as the original ones. When the marginal distributions of the
original vector and the approximating one are the same, there are some
results pertaining weak convergence metrized by integral probability metrics
that have not been explored before. The goal of this paper is to relate the
minimal generator of that metrizes weak convergence under an integral
probability metric to more general classes of test functions. As a few
examples show, this is particularly important when dealing with convergence
rates of moments between the original process and the approximation. To give
a simple illustration of the applications to be discussed below, consider $%
X_{n}:=\left( X_{1n},X_{2n}\right) $ converging weakly to $X:=\left(
X_{1},X_{2}\right) \in \mathbb{R}^{2}$. This implies that for any continuous
bounded function, say $f$, $\mathbb{E}f\left( X_{n}\right) \rightarrow 
\mathbb{E}f\left( X\right) .$ If $f\left( X\right) $ is only uniformly
integrable, the rate of convergence is not the one obtained under the Dudley
metric (bounded Lipschitz metric), as we need to truncate. However, if $%
law\left( X_{in}\right) =law\left( X_{i}\right) $ $\left( i=1,2\right) $,
under some conditions on $f$, the results of this paper show that truncation
is not necessary and a tighter bound can be achieved. In particular, $f$ may
even be discontinuous for the convergence of moments of $f$ to be satisfied
and the marginal laws can be arbitrary as long as they are in common. No
structure is imposed on $law\left( X_{n}\right) $ and $law\left( X\right) $
apart from the marginals being in common and the assumed weak convergence.

Below we give some background information also to set up the notation
(Section 2). Then, the result of the paper is stated (Section 3) followed by
some applications (Section 4). The proofs are in Sections 5 and 6.

\section{Background Material and Notation}

Suppose $\mathbb{P}\ $and $\mathbb{Q}$\ are two laws. For a class of
functions $\mathfrak{F}$ we define the integral metric 
\begin{equation}
d_{\mathfrak{F}}\left( \mathbb{P},\mathbb{Q}\right) :=\sup_{f\in \mathfrak{F}%
}\left\vert \dint f\left( x\right) d\left( \mathbb{P}-\mathbb{Q}\right)
\left( x\right) \right\vert .  \label{EQIntProbMetric}
\end{equation}%
These metrics have been used by Rachev (1991) to study the stability of
stochastic systems as well as by M\"{u}ller (1997). Suppose $f$\ is a
function with domain in the metric space $\left( \mathcal{X},r\right) $.\
Define the supremum norm $\left\Vert f\right\Vert _{\infty }:=\sup_{x\in 
\mathcal{X}}\left\vert f\left( x\right) \right\vert $ and the Lipschitz norm
as $\left\Vert f\right\Vert _{L}:=\sup_{x,y\in \mathcal{X}}\left\vert
f\left( x\right) -f\left( y\right) \right\vert /r\left( x,y\right) $. Define
the bounded Lipschitz metric of $f$\ as the semimetric $\left\Vert
f\right\Vert _{BL}:=\left\Vert f\right\Vert _{\infty }+\left\Vert
f\right\Vert _{L}.$ We say that a function $f$ belongs to the class $%
BL_{1}:=BL_{1}\left( \mathcal{X}\right) $ if $\left\Vert f\right\Vert
_{BL}\leq 1$. The metric space $\left( \mathcal{X},r\right) $\ is implicit
in this definition.

When $\mathfrak{F}:=BL_{1},$ and $\left( \mathbb{P}_{n}\right) _{n\in 
\mathbb{N}}$\ $\left( \mathbb{Q}_{n}\right) _{n\in \mathbb{N}}$ are two
sequences of laws defined on $\left( \mathcal{X},r\right) $, (\ref%
{EQIntProbMetric})\ can be used to metrize weak convergence of $\mathbb{P}%
_{n}$ and $\mathbb{Q}_{n}$ to each other (and often, $\mathbb{Q}_{n}=\mathbb{%
Q}$ is fixed and understood as the limiting value of $\left( \mathbb{P}%
_{n}\right) _{n\in \mathbb{N\cup }\left\{ \infty \right\} }$). Suppose $%
\mathcal{X}\subseteq \mathbb{R}^{K}$ and $\mathbb{P}_{n}$ and $\mathbb{Q}%
_{n} $ are the laws of $X_{n}:=\left( X_{n1},...,X_{nK}\right) $ and $%
X_{n}^{\prime }:=\left( X_{n1}^{\prime },...,X_{nK}^{\prime }\right) $\ such
that $law\left( X_{nk}\right) =law\left( X_{nk}^{\prime }\right) $ so that $%
\mathbb{P}_{n}$ and $\mathbb{Q}_{n}$ have the same univariate marginals.
This paper gives a result about the generator $\mathfrak{F}$ of (\ref%
{EQIntProbMetric}) when $\mathbb{P}_{n}$ and $\mathbb{Q}_{n}$ satisfy the
just mentioned restriction. The space $\mathcal{X}$ will be equipped with
the $l_{p}$ distance denoted by $r_{p}\left( x,y\right) :=\left(
\sum_{k=1}^{K}\left\vert x_{k}-y_{k}\right\vert ^{p}\right) ^{1/p}$ ($p\in %
\left[ 1,\infty \right] $, with the obvious modification for $p=\infty $)$.$%
\ Finally, $M_{1}:=M_{1}\left( \mathcal{X}\right) $ is the class of
coordinatewise nondecreasing functions from $\mathcal{X}$\ to $\left[ 0,1%
\right] .$\ 

\section{Statement of Result}

When the marginals between two laws are in common, in the sense defined
above, we have the following relationship between the generators $BL_{1}$
and $M_{1}.$

\begin{theorem}
\label{TheoremTopologicalEquivalence}Suppose $\left( \mathbb{P}_{n}\right)
_{n\in \mathbb{N}}$\ and $\left( \mathbb{Q}_{n}\right) _{n\in \mathbb{N}}$
have common marginals. Then, $d_{BL_{1}}\left( \mathbb{P}_{n},\mathbb{Q}%
_{n}\right) \rightarrow 0$ if and only if $d_{M_{1}}\left( \mathbb{P}_{n},%
\mathbb{Q}_{n}\right) \rightarrow 0$.
\end{theorem}

\begin{remark}
\label{RemarkUniformconvergence}Surprisingly, the two metrics $d_{BL_{1}}$
and $d_{M_{1}}$\ are topological equivalent when $\mathbb{P}_{n}$ and $%
\mathbb{Q}_{n}$ have same marginals. Since the indicator of $[x,\infty
):=[x_{1},\infty )\times \cdots \times \lbrack x_{K},\infty )$\ is
coordinatewise increasing, 
\begin{equation*}
\sup_{x\in \mathcal{X}}\left\vert \mathbb{P}_{n}\left( [x,\infty )\right) -%
\mathbb{Q}_{n}\left( [x,\infty )\right) \right\vert \leq d_{M_{1}}\left( 
\mathbb{P}_{n},\mathbb{Q}_{n}\right) ,
\end{equation*}%
so that weak convergence of sequences with common marginals is equivalent to
the above uniform convergence. By the Portmanteau Theorem, we know that this
is not true in general.
\end{remark}

In applications it is important to use functions that are not necessarily
continuous. However, for actual calculations, it is much easier to compute $%
d_{BL_{1}}$. We give a simple illustrative example.

\begin{example}
Suppose 
\begin{eqnarray*}
Y_{n} &=&\dsum_{s=0}^{\infty }a_{s}Z_{n-s} \\
Y_{n}^{\prime } &=&\dsum_{s=0}^{n-1}a_{s}Z_{n-s}+\dsum_{s=n}^{\infty
}a_{s}Z_{n-s}^{\prime },
\end{eqnarray*}%
where $\left( Z_{t}\right) _{t\in \mathbb{Z}}$ are iid integrable random
variables and $\left( Z_{t}^{\prime }\right) _{t\in \mathbb{Z}}$\ is an
independent copy of $\left( Z_{t}\right) _{t\in \mathbb{Z}}$. Hence $%
Y_{n}^{\prime }$ is independent of $Y_{0}$.\ Consider the laws of $%
X_{n}:=\left( Y_{0},Y_{n}\right) $\ and $X_{n}^{\prime }:=\left(
Y_{0},Y_{n}^{\prime }\right) ,$\ say $\mathbb{P}_{n}$ and $\mathbb{Q}_{n}$,
where $X_{n}$ and $X_{n}^{\prime }$\ take values in $\left( \mathbb{R}%
^{2},r_{1}\right) .$\ We have 
\begin{equation*}
d_{BL_{1}}\left( \mathbb{P}_{n},\mathbb{Q}_{n}\right) \leq \mathbb{E}%
r_{1}\left( X_{n},X_{n}^{\prime }\right) =\mathbb{E}\left\vert
Y_{n}-Y_{n}^{\prime }\right\vert \leq 2\mathbb{E}\left\vert Z_{0}\right\vert
\dsum_{s=n}^{\infty }\left\vert a_{s}\right\vert ,
\end{equation*}%
by the Lipschitz condition. If $\dsum_{s=n}^{\infty }\left\vert
a_{s}\right\vert \rightarrow 0$ as $n\rightarrow \infty ,$\ Theorem \ref%
{TheoremTopologicalEquivalence} gives $d_{M_{1}}\left( \mathbb{P}_{n},%
\mathbb{Q}_{n}\right) \rightarrow 0$\ as well. In words, this means that the
law of $X_{n}$ converges to the product of the marginals ($X_{n}^{\prime }$
has independent components) and the convergence holds for the expectation of
more general classes of functions than $BL_{1}$ (e.g. Remark \ref%
{RemarkUniformconvergence}).
\end{example}

As the previous example shows, it is also important to relate the two
metrics in order to derive rates of convergence.

\begin{theorem}
\label{TheoremRelationship}Suppose $\left( \mathbb{P}_{n}\right) _{n\in 
\mathbb{N}}$\ and $\left( \mathbb{Q}_{n}\right) _{n\in \mathbb{N}}$ are as
in Theorem \ref{TheoremTopologicalEquivalence} and recall that $\mathbb{P}%
_{n}$ and $\mathbb{Q}_{n}$ are defined on $\left( \mathcal{X\subseteq }%
\mathbb{R}^{K},r_{p}\right) $ .\ Then, 
\begin{equation*}
d_{M_{1}}\left( \mathbb{P}_{n},\mathbb{Q}_{n}\right) \leq \min \left\{
2^{3/2}K^{\left( p-1\right) /\left( 2p\right) }d_{BL_{1}}\left( \mathbb{P}%
_{n},\mathbb{Q}_{n}\right) ^{1/2},1\right\} .
\end{equation*}
\end{theorem}

\section{Two Applications}

We give two applications. In both cases, we show that relating the integral
probability metric of the class $BL_{1}$ to the class $M_{1}$\ allows us to
derive upper bounds for integral probability metrics based on unbounded
classes of functions, not necessarily continuous. In particular we shall
define $MRC_{+}:=MRC_{+}\left( \mathbb{R}\right) $\ to be the class of
positive nondecreasing functions on $\mathbb{R}$, continuous from the right.

\subsection{Integral Probability Metrics for Products of Unbounded Positive
Functions}

Suppose 
\begin{equation*}
\mathfrak{G}:=\left\{ \dprod_{k=1}^{K}g_{k}:g_{k}\in
MRC_{+},k=1,...,K\right\}
\end{equation*}%
\ is a class of functions such that each element $g$ is the pointwise
product of positive nondecreasing functions continuous from the right. Even
if $g\in \mathfrak{G}$ is unbounded, the integral probability metric with
respect to $\mathfrak{G}$ can be bounded in terms of $M_{1}\left( \mathbb{R}%
^{K}\right) $ and $BL_{1}\left( \mathbb{R}^{K}\right) $ by Theorem \ref%
{TheoremRelationship}.

\begin{corollary}
\label{CorollaryUnboundedClass}Suppose $\mathbb{P}$, $\mathbb{Q}$\ are the
laws of $X$, $X^{\prime }\in \mathbb{R}^{K}$ having same univariate
marginals. Set $\theta :=d_{M_{1}\left( \mathbb{R}^{K}\right) }\left( 
\mathbb{P},\mathbb{Q}\right) /2$, then%
\begin{equation*}
d_{\mathfrak{G}}\left( \mathbb{P},\mathbb{Q}\right) \leq \dint_{0}^{\theta
}\dprod_{k=1}^{K}Q_{g_{k}}\left( s\right) ds,
\end{equation*}%
where $Q_{g_{k}}\left( s\right) =\inf \left\{ x:\Pr \left( g_{k}\left(
X_{k}\right) >x\right) \leq s\right\} $.
\end{corollary}

\begin{remark}
By simple approximating arguments, we can extend Corollary \ref%
{CorollaryUnboundedClass} from $\mathfrak{G}$ to the convex hull of $%
\mathfrak{G}\cup \left( -\mathfrak{G}\right) ,$ where $g\in \left( -%
\mathfrak{G}\right) $\ if $-g\in \mathfrak{G}$. Details are left to the
reader.
\end{remark}

\subsection{Covariance Inequalities}

Suppose, $Y$ and $Z$ are two random variables. Define 
\begin{equation}
\alpha :=\alpha \left( Y,Z\right) :=2\sup_{z,y\in \mathbb{R}}\left\vert \Pr
\left( Y>y,Z>z\right) -\Pr \left( Y>y\right) \Pr \left( Z>z\right)
\right\vert ,  \label{EQSimpleStrongMixing}
\end{equation}%
which is a simple version of the strong mixing coefficient 
\begin{equation}
\alpha \left( \mathcal{A},\mathcal{B}\right) :=2\sup \left\{ Cov\left(
I_{A},I_{B}\right) :A\in \mathcal{A},B\in \mathcal{B}\right\} ,
\label{EQStrongMixingDef}
\end{equation}%
where $\mathcal{A}$ and $\mathcal{B}$ are two sigma algebras. Strong mixing
is a commonly used dependence condition and limit theorem for dependent
random variables are often proved using it. The coefficient $\alpha $ is the
basis of Rio's covariance inequality (e.g. Rio 2000, Theorem 1.1), which is
known to be sharp. Using the results of this paper we can prove the
following extension of the just mentioned inequality.

\begin{corollary}
\label{CorollaryRio}Suppose $\left( Y,Z\right) $ is a random vector with
values in $\left( \mathbb{R}^{2},r_{1}\right) $. Set 
\begin{equation*}
\theta :=2\min \left\{ \left\vert 8d_{BL_{1}}\left( \mathbb{P}_{\left(
Y,Z\right) },\mathbb{P}_{Y}\mathbb{P}_{Z}\right) \right\vert
^{1/2},1\right\} ,
\end{equation*}%
where $\mathbb{P}_{\left( Y,Z\right) }$ is the law of $\left( Y,Z\right) $\
and $\mathbb{P}_{Y}\mathbb{P}_{Z}$\ it the product of the marginals. Suppose 
$g_{Y},g_{Z}\in MRC_{+}$.\ Then, 
\begin{equation*}
Cov\left( g_{Y}\left( Y\right) ,g_{Z}\left( Z\right) \right) \leq
2\dint_{0}^{\theta }Q_{g_{Y}}\left( s\right) Q_{g_{Z}}\left( s\right) ds,
\end{equation*}%
where $Q_{g_{Y}}\left( s\right) :=\inf \left\{ y:\Pr \left( \left\vert
g_{Y}\left( Y\right) \right\vert >y\right) \leq s\right\} $ and $%
Q_{g_{Z}}\left( s\right) :=\inf \left\{ z:\Pr \left( \left\vert g_{Z}\left(
Z\right) \right\vert >z\right) \leq s\right\} $.
\end{corollary}

\begin{remark}
To show a similar result, Doukhan and Louhichi (1999, Lemma 1) require the
marginals to be Lipschitz continuous. Theorem \ref{TheoremRelationship}
shows that Lipschitz continuity is not necessary. In general the strong
mixing condition is satisfied if $d_{BL_{1}}\rightarrow 0$, and we restrict $%
A$ and $B$ in (\ref{EQStrongMixingDef}) to be half open intervals, as in (%
\ref{EQSimpleStrongMixing}).
\end{remark}

\section{Weak Convergence with Common Marginals}

Suppose $X:=\left( X_{1},...,X_{K}\right) \in \mathcal{X}\subseteq \mathbb{R}%
^{K}$\ is a vector of random variables with law $\mathbb{P}$\ and $%
P_{k}:=law\left( X_{k}\right) $ $\left( k=1,...,K\right) .$\ Then, there
exists a function $C:\left[ 0,1\right] ^{K}\rightarrow \left[ 0,1\right] $
such that 
\begin{equation}
\Pr \left( X_{1}\in A_{1},...,X_{K}\in A_{K}\right) =C\left( P_{1}\left(
A_{1}\right) ,...,P_{K}\left( A_{K}\right) \right)
\label{EQCopulaRepresentation}
\end{equation}%
for $A_{k}=\left[ -\infty ,x_{k}\right] \subset \mathbb{R},$ $\left(
k=1,...,K\right) $. The function $C$ is called copula function and the above
representation holds even for sets more general than $A_{k}$\ (Scarsini,
1989). Suppose $C^{\prime }$ is a copula distinct from $C$ above. The
distance between these two copulae can be measured in terms of an integral
probability metric generated by a class of functions $\mathfrak{F}$, i.e. $%
d_{\mathfrak{F}}\left( C,C^{\prime }\right) .$ Define%
\begin{equation*}
\tilde{F}_{k}\left( x,v\right) :=\Pr \left( X_{k}<x\right) +v\Pr \left(
X_{k}=x\right) ,\text{ }v\in \left[ 0,1\right] \text{ and }x\in \mathbb{R}.
\end{equation*}%
Define the operator $\mathcal{Q}^{\ast }:\mathbb{R}^{K}\rightarrow \left[ 0,1%
\right] ^{K}$ such that 
\begin{equation*}
\mathcal{Q}^{\ast }X:=\left( \tilde{F}_{1}\left( X_{1},V_{1}\right) ,...,%
\tilde{F}_{K}\left( X_{K},V_{K}\right) \right) ,
\end{equation*}%
where $\left( V_{1},...,V_{K}\right) $\ are iid $\left[ 0,1\right] $ uniform
random variables independent of $X$. Then, $U:=\mathcal{Q}^{\ast }X$ is a
random vector with uniform $\left[ 0,1\right] $ marginals (R\"{u}schendorf
and de Valk, 1993, Proposition 1).\ If $P_{k}$\ $\left( k=1,...,K\right) $\
is absolutely continuous with respect to the Lebesgue measure,\ from (\ref%
{EQCopulaRepresentation})\ we deduce $U$ has law $C.$ Define also the
operator $\mathcal{Q}:\left[ 0,1\right] ^{K}\rightarrow \mathbb{R}^{K}$ such
that for some $\left[ 0,1\right] ^{K}$\ uniform random vector $U$ with law $%
C $,%
\begin{equation*}
\mathcal{Q}U:=\left( P_{1}^{-1}\left( U_{1}\right) ,...,P_{K}^{-1}\left(
U_{K}\right) \right) ,\text{ (where }P_{k}^{-1}\left( u\right) :=\inf
\left\{ x:\Pr \left( X\leq x\right) \geq u\right\} \text{)}
\end{equation*}%
so that $\mathcal{Q}U\overset{d}{=}X$ ($\overset{d}{=}$ is equality in
distribution)$.$ Hence, if the marginals of $\mathcal{Q}U$ are not
continuous, the copula is not unique (e.g. Scarsini, 1989). Here, the term
copula will refer to the law of $U\overset{d}{=}\mathcal{Q}^{\ast }X$, which
is unique.

Denote $M_{1}\left( \left[ 0,1\right] ^{K}\right) $ to be the class of
functions $M_{1}$, but with support in $\left[ 0,1\right] ^{K}$; similarly,
define $BL_{1}\left( \left[ 0,1\right] ^{K}\right) $. Moreover define 
\begin{equation*}
\mathfrak{B}:=\left\{ g:\left[ 0,1\right] ^{K}\rightarrow \mathbb{R}%
:g:=f\circ \mathcal{Q},f\in BL_{1}\left( \mathcal{X}\right) \right\}
\end{equation*}%
so that $\mathfrak{B}$ is the class of functions obtained from functions in $%
BL_{1}\left( \mathcal{X}\right) $\ using the composition of $f$ with the
transformation $\mathcal{Q}$.

\begin{lemma}
\label{LemmaEquivalenceProbabilityCopulaIntegral}Suppose $X\ $and $X^{\prime
}$ are random variables with values in $\mathcal{X}\subseteq \mathbb{R}^{K}$
and laws $\mathbb{P}$ and $\mathbb{Q}$\ having common marginals and with
copula $C$ and $C^{\prime }$ respectively. Then 
\begin{equation*}
d_{M_{1}\left( \mathcal{X}\right) }\left( \mathbb{P},\mathbb{Q}\right)
=d_{M_{1}\left( \left[ 0,1\right] ^{K}\right) }\left( C,C^{\prime }\right)
\end{equation*}%
and 
\begin{equation*}
d_{BL_{1}\left( \mathcal{X}\right) }\left( \mathbb{P},\mathbb{Q}\right) =d_{%
\mathfrak{B}}\left( C,C^{\prime }\right) \geq d_{BL_{1}\left( \left[ 0,1%
\right] ^{K}\right) }\left( C,C^{\prime }\right) .
\end{equation*}
\end{lemma}

\begin{proof}
Note that if $f\in M_{1}\left( \mathcal{X}\right) $ then $f\circ \mathcal{Q}%
\in M_{1}\left( \left[ 0,1\right] ^{K}\right) $ and if $g\in M_{1}\left( %
\left[ 0,1\right] ^{K}\right) $ then $g\circ \mathcal{Q}^{-1}\in M_{1}\left( 
\mathcal{X}\right) ,$ where $\mathcal{Q}^{-1}X:=\left( P_{1}\left(
X_{1}\right) ,...,P_{K}\left( X_{K}\right) \right) $\ and $P_{k}\left(
x\right) :=P_{k}\left( \left[ -\infty ,x\right] \right) $. This follows
because $\mathcal{Q}$\ is nondecreasing. Since $\mathcal{Q}$\ is a
measurable map (if we consider the sigma algebras generated by Borel sets of 
$\mathcal{X}$ and $\left[ 0,1\right] ^{K}$)\ by a change of variables,%
\begin{equation*}
d_{M_{1}\left( \mathcal{X}\right) }\left( \mathbb{P},\mathbb{Q}\right)
=\sup_{f\in M_{1}\left( \mathcal{X}\right) }\left\vert \dint fd\left( 
\mathbb{P}-\mathbb{Q}\right) \right\vert =\sup_{f\in M_{1}\left( \mathcal{X}%
\right) }\left\vert \dint f\circ \mathcal{Q}d\left( C-C^{\prime }\right)
\right\vert =d_{M_{1}\left( \left[ 0,1\right] ^{K}\right) }\left(
C,C^{\prime }\right) .
\end{equation*}%
The last equality follows because $f=g\circ \mathcal{Q}^{-1}\in M_{1}\left( 
\mathcal{X}\right) $ (for some $g\in M_{1}\left( \left[ 0,1\right]
^{K}\right) $), as mentioned above. The second part of the lemma follows by
change of variables. We only need to note that functions in $\mathfrak{B}$
are bounded, but not necessarily Lipschitz unless the marginals of $X$ are
Lipschitz continuos. Hence, $BL_{1}\left( \left[ 0,1\right] ^{K}\right)
\subseteq \mathfrak{B}$.
\end{proof}

For $U\in \left[ 0,1\right] ^{K}$ define $\bar{U}=\left(
1-U_{1},...,1-U_{K}\right) $ and similarly for $u\in \left[ 0,1\right] ^{K}$
define $\bar{u}$. Then,$\ \hat{C}\left( u\right) :=\Pr \left( \bar{U}\leq
u\right) $\ is the copula of $\bar{U}$ so that $\bar{C}\left( u\right) :=%
\hat{C}\left( \bar{u}\right) =\Pr \left( U\geq u\right) $\ is the survival
copula (vector inequalities are meant to hold elementwise). We recall that
the copula is Lipschitz of order and constant one with respect to the $l_{1}$
norm, 
\begin{equation}
\left\vert C\left( u_{1},...,u_{K}\right) -C\left( v_{1},...,v_{K}\right)
\right\vert \leq r_{1}\left( u,v\right) =\dsum_{k=1}^{K}\left\vert
u_{k}-v_{k}\right\vert .  \label{EQCopulaLipschitzContinuity1}
\end{equation}%
Since $\hat{C}\left( u\right) $ is a copula, (\ref%
{EQCopulaLipschitzContinuity1}) also implies%
\begin{equation}
\left\vert \bar{C}\left( u\right) -\bar{C}\left( v\right) \right\vert
=\left\vert \hat{C}\left( \bar{u}\right) -\hat{C}\left( \bar{v}\right)
\right\vert \leq r_{1}\left( \bar{u},\bar{v}\right) =r_{1}\left( u,v\right) .
\label{EQCopulaLipschitzContinuity2}
\end{equation}%
\ Hence, we have the following.

\begin{lemma}
\label{LemmaMaximalGenerator}%
\begin{equation*}
\sup_{f\in M_{1}\left( \left[ 0,1\right] ^{K}\right) }\left\vert \int
fd\left( C-C^{\prime }\right) \right\vert =\sup_{u\in \left[ 0,1\right]
^{K}}\left\vert \bar{C}\left( u\right) -\bar{C}^{\prime }\left( u\right)
\right\vert .
\end{equation*}
\end{lemma}

\begin{proof}
For simplicity, suppose $f\left( 0_{K}\right) =0$ ($0_{K}$ is the $K$
dimensional vector of zeros)$.$ Define 
\begin{equation*}
f_{m}:=\sum_{i\in \mathcal{I}_{m}}\Delta _{m}f_{i}I\left\{ u\in \left[
u_{i},1\right] \right\} ,
\end{equation*}%
where $\left[ u_{i},1\right] :=\left[ u_{i_{1}},1\right] \times \cdots
\times \left[ u_{i_{K}},1\right] \subset \left[ 0,1\right] ^{K},$ $%
u_{i_{k}}\leq u_{i_{k}+1}$,\ and $\Delta _{m}f_{i}\geq 0$ such that $%
\sum_{i\in \mathcal{I}_{m}}\Delta _{m}f_{i}\leq 1,$ $i\in \mathcal{I}%
_{m}\subset \mathbb{N}^{K}$\ and $\#\mathcal{I}_{m}=2^{m}$ ($I\left\{
...\right\} $ is the indicator function)$.$\ Choosing $\left( u_{i}\right)
_{i\in \mathcal{I}_{m}}$ such that $\Delta _{m}f_{i}\leq c2^{-m}$ (for some
bounded absolute constant $c$), $\mathbb{E}\left\vert f_{m}-f\right\vert
\rightarrow 0$ as $m\rightarrow \infty $ for any $f\in M_{1}$. Such a
construction and its convergence can be shown by simple manipulation of the
usual approximation of functions in terms of indicator functions (e.g. see
the proof of Theorem 13.5 in Billingsley, 1995) together with monotonicity.
Then,%
\begin{eqnarray*}
\left\vert \int f_{m}d\left( C-C^{\prime }\right) \right\vert  &=&\left\vert
\sum_{i\in \mathcal{I}_{m}}\Delta _{m}f_{i}\int I\left\{ u\in \left[ u_{i},1%
\right] \right\} d\left( C-C^{\prime }\right) \left( u\right) \right\vert  \\
&\leq &\sum_{i\in \mathcal{I}_{m}}\Delta _{m}f_{i}\max_{i\in \mathcal{I}%
_{m}}\left\vert \int I\left\{ u\in \left[ u_{i},1\right] \right\} d\left(
C-C^{\prime }\right) \left( u\right) \right\vert .
\end{eqnarray*}%
By continuity of $C$ (i.e. (\ref{EQCopulaLipschitzContinuity1})),\ the above
display is unaffected if $\left[ u_{i},1\right] $ in the definition of $f_{m}
$ does not include the left boundary in any/some of the coordinates. Hence,
we can regard $\sup_{m}f_{m}$ to be a nondecreasing function (not
necessarily right continuous). Hence, taking $\sup $ with respect to $m$ and 
$f$, we deduce 
\begin{equation*}
\sup_{f\in M_{1}}\left\vert \int fd\left( C-C^{\prime }\right) \right\vert
\leq \sup_{u\in \left[ 0,1\right] ^{K}}\left\vert \Pr \left( U\geq u\right)
-\Pr \left( U^{\prime }\geq u\right) \right\vert =\sup_{u\in \left[ 0,1%
\right] ^{K}}\left\vert \bar{C}\left( u\right) -\bar{C}^{\prime }\left(
u\right) \right\vert ,
\end{equation*}%
where $U$ and$\ U^{\prime }$\ are uniform $\left[ 0,1\right] ^{K}$\ random
vectors with law $C$ and $C^{\prime }$ respectively. The lower bound follows
noting that the indicator of $\left[ u_{i},1\right] \subset \left[ 0,1\right]
^{K}$ is coordinatewise increasing, where the same remark about the
boundaries applies.
\end{proof}

\begin{proof}[Proof of Theorem \protect\ref{TheoremTopologicalEquivalence}]
Suppose $d_{BL_{1}\left( \mathcal{X}\right) }\left( \mathbb{P}_{n},\mathbb{Q}%
_{n}\right) \rightarrow 0.$\ By Lemma \ref%
{LemmaEquivalenceProbabilityCopulaIntegral},\ $d_{BL_{1}\left( \left[ 0,1%
\right] ^{K}\right) }\left( C_{n},C_{n}^{\prime }\right) \rightarrow 0$\
where $C_{n}$ and $C_{n}^{\prime }$\ are the copulae corresponding to $%
\mathbb{P}_{n}$ and $\mathbb{Q}_{n}.$ Then, by (\ref%
{EQCopulaLipschitzContinuity2}) (i.e. continuity), 
\begin{equation}
\sup_{u\in \left[ 0,1\right] ^{K}}\left\vert \bar{C}_{n}\left( u\right) -%
\bar{C}_{n}^{\prime }\left( u\right) \right\vert \rightarrow 0
\label{EQSupConvergence}
\end{equation}%
and $d_{M_{1}\left( \left[ 0,1\right] ^{K}\right) }\left(
C_{n},C_{n}^{\prime }\right) \rightarrow 0$\ by Lemma \ref%
{LemmaMaximalGenerator}. Finally, Lemma \ref%
{LemmaEquivalenceProbabilityCopulaIntegral}\ implies $d_{M_{1}\left( 
\mathcal{X}\right) }\left( \mathbb{P}_{n},\mathbb{Q}_{n}\right) \rightarrow
0.$

Now suppose $d_{M_{1}\left( \mathcal{X}\right) }\left( \mathbb{P}_{n},%
\mathbb{Q}_{n}\right) \rightarrow 0.\ $Since the class of indicators of
right half spaces are in $M_{1}\left( \mathcal{X}\right) $, 
\begin{equation*}
\sup_{x\in \mathcal{X}}\left\vert \Pr \left( X_{n}\geq x\right) -\Pr \left(
X_{n}^{\prime }\geq x\right) \right\vert \leq d_{M_{1}\left( \mathcal{X}%
\right) }\left( \mathbb{P},\mathbb{Q}\right) \rightarrow 0,
\end{equation*}%
and we must have $d_{BL_{1}\left( \mathcal{X}\right) }\left( \mathbb{P}_{n},%
\mathbb{Q}_{n}\right) \rightarrow 0$, by weak convergence (because in $%
\mathcal{X}\subseteq \mathbb{R}^{K}$ the sets $\left[ x,\infty \right] :=%
\left[ x_{1},\infty \right] \times \cdots \times \left[ x_{K},\infty \right]
\subset \mathbb{R}^{K}$ are a separating class).
\end{proof}

\begin{proof}[Proof of Theorem \protect\ref{TheoremRelationship}]
At first we show that the bound holds for $C_{n}$ and $C_{n}^{\prime }.$\
Let $f\left( u\right) :=I\left\{ u\in \left[ v,1\right] \right\} $ and $%
g_{\epsilon }\left( u\right) :=\left[ 1-r_{p}\left( \left[ v,1\right]
,u\right) /\epsilon \right] \vee 0$, $u\in \left[ 0,1\right] ^{K}$ ($\left[
v,1\right] :=\left[ v_{1},1\right] \times \cdots \times \left[ v_{K},1\right]
\subset \left[ 0,1\right] ^{K}$). Hence, $g_{\epsilon }$\ is a Lipschitz
approximation of $f$ (as $\epsilon \rightarrow 0$) and such that $\left\Vert
g_{\epsilon }\right\Vert _{BL}\leq \left( 1+\epsilon ^{-1}\right) $. Noting
that $g_{\epsilon }\geq f,$\ consider the following bound, where dependence
on the argument $u$ is suppressed,%
\begin{eqnarray}
\int fd\left( C_{n}-C_{n}^{\prime }\right) &\leq &\int g_{\epsilon
}dC_{n}-\int fdC_{n}^{\prime }=\int g_{\epsilon }d\left( C_{n}-C_{n}^{\prime
}\right) +\int \left( g_{\epsilon }-f\right) dC_{n}^{\prime }  \notag \\
&\leq &\left( 1+\epsilon ^{-1}\right) d_{BL_{1}\left( \left[ 0,1\right]
^{K}\right) }\left( C_{n},C_{n}^{\prime }\right) +\int \left( g_{\epsilon
}-f\right) dC_{n}^{\prime }  \notag \\
&&\text{[by the previous remarks about }g_{\epsilon }\text{]}  \notag \\
&\leq &\left( 1+\epsilon ^{-1}\right) d_{BL_{1}\left( \left[ 0,1\right]
^{K}\right) }\left( C_{n},C_{n}^{\prime }\right) +r_{1}\left( u^{\prime
},v\right) ,  \label{EQMBLBound}
\end{eqnarray}%
for some $u^{\prime }\in \left[ 0,1\right] ^{K}$\ using (\ref%
{EQCopulaLipschitzContinuity2}).\ By construction, $u^{\prime }$ and $v$ are
such that%
\begin{equation*}
r_{p}\left( v,u^{\prime }\right) =\left( \dsum_{k=1}^{K}\left\vert
u_{k}^{\prime }-v_{k}\right\vert ^{p}\right) ^{1/p}\leq \epsilon ,
\end{equation*}%
implying, by the relations between means (i.e. $r_{1}\left( u,v\right)
K^{-1}\leq r_{p}\left( u,v\right) K^{-1/p}$), 
\begin{equation*}
r_{1}\left( u^{\prime },v\right) \leq K^{1-1/p}r_{p}\left( u^{\prime
},v\right) \leq K^{\left( p-1\right) /p}\epsilon .
\end{equation*}%
\ Substitute this bound in (\ref{EQMBLBound}). Noting that $\epsilon \leq 1,$%
\ choosing $\epsilon $\ to equate the two terms gives 
\begin{equation*}
d_{M_{1}\left( \left[ 0,1\right] ^{K}\right) }\left( C_{n},C_{n}^{\prime
}\right) \leq 2^{3/2}K^{\left( p-1\right) /\left( 2p\right) }d_{BL_{1}\left( %
\left[ 0,1\right] ^{K}\right) }\left( C_{n},C_{n}^{\prime }\right) ^{1/2},
\end{equation*}%
where we have taken $\sup $ with respect to $f$ on both sides ($f$ in the
class of indicators of right half open intervals $\left[ v,1\right] \subset %
\left[ 0,1\right] ^{K}$ as defined above) and we have used Lemma \ref%
{LemmaMaximalGenerator}. Then, by Lemma \ref%
{LemmaEquivalenceProbabilityCopulaIntegral}, the previous display implies 
\begin{eqnarray*}
d_{M_{1}\left( \mathcal{X}\right) }\left( \mathbb{P}_{n},\mathbb{Q}%
_{n}\right) &=&d_{M_{1}\left( \left[ 0,1\right] ^{K}\right) }\left(
C_{n},C_{n}^{\prime }\right) \leq 2^{3/2}K^{\left( p-1\right) /\left(
2p\right) }d_{BL_{1}\left( \left[ 0,1\right] ^{K}\right) }\left(
C_{n},C_{n}^{\prime }\right) ^{1/2} \\
&\leq &2^{3/2}K^{\left( p-1\right) /\left( 2p\right) }d_{BL_{1}\left( 
\mathcal{X}\right) }\left( \mathbb{P},\mathbb{Q}\right) ^{1/2}.
\end{eqnarray*}%
By Lemma \ref{LemmaMaximalGenerator}, we deduce $d_{M_{1}\left( \left[ 0,1%
\right] ^{K}\right) }\left( C_{n},C_{n}^{\prime }\right) \leq 1$ and the
result follows.
\end{proof}

\section{Proof of Corollaries}

\begin{proof}[Proof of Corollary \protect\ref{CorollaryUnboundedClass}]
For any positive function $g_{k}:\mathbb{R}\rightarrow \mathbb{R}$%
\begin{equation*}
g_{k}\left( x\right) =\dint_{0}^{\infty }I\left\{ g_{k}\left( x\right)
>s\right\} ds.
\end{equation*}%
Hence, suppressing the argument $x$ for ease of notation,%
\begin{eqnarray*}
\text{I} &:&=\dint g_{1}\cdots g_{K}d\left( \mathbb{P}-\mathbb{Q}\right) \\
&=&\dint \dprod_{k=1}^{K}\dint_{0}^{\infty }I\left\{ g_{k}>s_{k}\right\}
ds_{k}d\left( \mathbb{P}-\mathbb{Q}\right) \\
&=&\int_{\mathbb{R}^{K}}ds_{1}\cdots ds_{K}\dint \dprod_{k=1}^{K}I\left\{
g_{k}>s_{k}\right\} d\left( \mathbb{P}-\mathbb{Q}\right)
\end{eqnarray*}%
by Fubini's Theorem. Then, 
\begin{eqnarray*}
\text{II} &:&=\dint \dprod_{k=1}^{K}I\left\{ g_{k}>s_{k}\right\} d\left( 
\mathbb{P}-\mathbb{Q}\right) \\
&\mathbb{\leq }&\left\vert \Pr \left( g_{k}\left( X_{k}\right)
>s_{k},k=1,..,K\right) -\Pr \left( g_{k}\left( X_{k}^{\prime }\right)
>s_{k},k=1,..,K\right) \right\vert \\
&=&\left\vert \Pr \left( X_{k}>g_{k}^{-1}\left( s_{k}\right)
,k=1,..,K\right) -\Pr \left( X_{k}^{\prime }>g_{k}^{-1}\left( s_{k}\right)
,k=1,..,K\right) \right\vert \\
&&\text{\lbrack because }g_{k}\text{ is nondecreasing and right continuous]}
\\
&\leq &d_{M_{1}\left( \mathbb{R}^{K}\right) }\left( \mathbb{P},\mathbb{Q}%
\right) =2\theta ,
\end{eqnarray*}%
and 
\begin{equation*}
\dint \dprod_{k=1}^{K}I\left\{ g_{k}>s_{k}\right\} d\left( \mathbb{P}-%
\mathbb{Q}\right) \mathbb{\leq }2\min \left( \Pr \left( g_{k}\left(
X_{k}\right) >s_{1}\right) ,...,\Pr \left( g_{K}\left( X_{K}\right)
>s_{K}\right) \right)
\end{equation*}%
because of the following reasons: the marginals of $\mathbb{P}$ and $\mathbb{%
Q}$\ are the same, and 
\begin{eqnarray*}
\dint \dprod_{k=1}^{K}I\left\{ g_{k}>s_{k}\right\} d\mathbb{P} &\mathbb{=}%
&\dint I\left\{ g_{l}>s_{l}\right\} \dprod_{\substack{ k=1  \\ k\neq l}}%
^{K}I\left\{ g_{k}>s_{k}\right\} d\mathbb{P} \\
&\leq &\min \left( \dint I\left\{ g_{1}>s_{1}\right\} d\mathbb{P},...,\dint
I\left\{ g_{K}>s_{K}\right\} d\mathbb{P}\right) \\
&&\text{[because the terms in the product take values in }\left[ 0,1\right] 
\text{]} \\
&=&\min \left( \Pr \left( g_{1}\left( X_{1}\right) >s_{1}\right) ,...,\Pr
\left( g_{K}\left( X_{K}\right) >s_{K}\right) \right) .
\end{eqnarray*}%
From the previous displays, 
\begin{eqnarray*}
\text{II} &\leq &2\min \left( \theta ,\Pr \left( g_{1}\left( X_{1}\right)
>s_{1}\right) ,...,\Pr \left( g_{K}\left( X_{K}\right) >s_{K}\right) \right)
\\
&=&2\int_{0}^{\theta }\dprod_{k=1}^{K}I\left\{ s<\Pr \left( g_{k}\left(
X_{k}\right) >s_{k}\right) \right\} ds=2\int_{0}^{\theta
}\dprod_{k=1}^{K}I\left\{ Q_{g_{k}}\left( s\right) >s_{k}\right\} ds,
\end{eqnarray*}%
which implies%
\begin{equation*}
\text{I}\leq \int_{\mathbb{R}^{K}}\text{II}ds_{1}\cdots ds_{K}\leq
\dint_{0}^{\theta }\dprod_{k=1}^{K}Q_{g_{k}}\left( s\right) ds,
\end{equation*}%
using Fubini's Theorem.
\end{proof}

\begin{proof}[Proof of Corollary \protect\ref{CorollaryRio}]
By Rio's covariance inequality,%
\begin{equation*}
Cov\left( g_{Y}\left( Y\right) ,g_{Z}\left( Z\right) \right) \leq
2\dint_{0}^{\alpha \left( g_{Y},g_{Z}\right) }Q_{g_{Y}}\left( s\right)
Q_{g_{Z}}\left( s\right) ds.
\end{equation*}%
Since $g_{Y}$ and $g_{Z}$ are nondecreasing and right continuous,\ 
\begin{eqnarray*}
&&\sup_{z,y\in \mathbb{R}}\left\vert \Pr \left( g_{Y}\left( Y\right)
>y,g_{Z}\left( Z\right) >z\right) -\Pr \left( g_{Y}\left( Y\right) >y\right)
\Pr \left( g_{Z}\left( Z\right) >z\right) \right\vert \\
&=&\sup_{z,y\in \mathbb{R}}\left\vert \Pr \left( Y>g_{Y}^{-1}\left( y\right)
,Z>g_{Z}^{-1}\left( z\right) \right) -\Pr \left( Y>g_{Y}^{-1}\left( y\right)
\right) \Pr \left( Z>g_{Z}^{-1}\left( z\right) \right) \right\vert ,
\end{eqnarray*}%
so, by Theorem \ref{TheoremRelationship}, 
\begin{equation*}
\alpha \left( g_{Y},g_{Z}\right) =\alpha \left( Y,Z\right) \leq
2d_{M_{1}}\left( \mathbb{P}_{\left( Y,Z\right) },\mathbb{P}_{Y}\mathbb{P}%
_{X}\right) \leq 2\min \left\{ \left\vert 8d_{BL_{1}}\left( \mathbb{P}%
_{\left( Y,Z\right) },\mathbb{P}_{Y}\mathbb{P}_{Z}\right) \right\vert
^{1/2},1\right\} .
\end{equation*}
\end{proof}

\end{document}